\documentclass[10pt]{article}
\usepackage{amsfonts}
\usepackage{authblk} 
\usepackage{hyperref} 
\usepackage{fancyhdr} 
\usepackage{amsthm} 
\usepackage{graphicx}
\usepackage{pstricks-add}
\usepackage{enumerate}
\usepackage{todonotes}
\usepackage{cancel}
\usepackage{float} 
\usepackage{tikz}
\usetikzlibrary{shapes,arrows}
\usepackage{amsmath}
\usepackage{subfig} 
\newtheorem{theorem}{Theorem}[section]

\newtheorem{lemma}[theorem]{Lemma}

\newcommand{\R}{\mathbb{R}}
\newcommand{\N}{\mathbb{N}}

\makeatletter
\@addtoreset{equation}{section}
\makeatother

\begin{document}

\title{A stochastic differential equation approach for an SIS model with non-linear incidence rate}

\author[1]{J.S Builes}
\affil[1]{Instituto de Matemáticas, Universidad de Antioquia, Medellín, Colombia}

\author[2]{Cristian F. Coletti}
\affil[2]{Universidade Federal do ABC, São Paulo, Brazil}

\author[3]{Leon A. Valencia}
\affil[3]{Instituto de Matemáticas, Universidad de Antioquia, Medellín, Colombia}

\maketitle

\begin{abstract}
In this paper, we study an analytically tractable SIS model with a non-linear incidence rate for the number of infectious individuals described through a stochastic differential equation (SDE). We guarantee the existence of a positive solution, and we study its regularity. We study the persistence and extinction regimes, and we give sufficient conditions under which the disease-free equilibrium point is an asymptotically stable equilibrium point with probability one. We provide sufficient conditions under which the model admits a unique stationary measure. Finally, we illustrate our findings using simulations.
\end{abstract}

\textbf{Keywords:} Stochastic differential equations,  S.I.S model, Extinction, Persistence, Stationary measure.
\section{Introduction}
Differential equations have many interesting applications in different fields of science. For instance, we use systems of differential equations and stochastic differential equations in the study of infections and disease transmission. We refer to these as deterministic epidemiological models and stochastic epidemiological models, respectively. One specific model, called the S.I.R. model, represents a population with susceptible, infected, and recovered individuals. The first mathematical model proposed for epidemics was given by Kermack-McKendrick \cite{Kermack}. Later, its generalizations include natality, mortality, loss of immunity, and some psychological effects on the population.
On the other hand, if the recovered population does not develop immunity to the disease, meaning all recovered individuals will be immediately susceptible to the disease again, the epidemiological model becomes an S.I.S. model. In \cite{Gray} and \cite{Xu}, a stochastic study of the S.I.S. model with births and deaths is presented. 

Now, consider the following system of deterministic differential equations

\begin{eqnarray}\label{sire1}
    dy &=&[-\beta xy+\gamma x +\mu(N- y)]dt,\nonumber\\
dx&=&[\beta xy -(\gamma+\mu)x]dt, \quad\text{ for $t\geq 0$}.
\end{eqnarray}
where $y$ is the density  of susceptible individuals and $x$ is the density of infected individuals. $N$ is the total population in terms of biomass, $\beta$ is the transmission rate of the disease, $\mu$ is the per-capital rate of births and deaths and $\gamma$ is the per-capital rate of recovery. In this model, we assume
 $x(t)+y(t)\equiv N$ and $\beta=\displaystyle\frac{\lambda}{N}$ and $\lambda$ is a per-capital rate of disease transmission. $\beta xy$ is called the incidence rate. In   
 \cite{MR529097} the authors modify the Kermac-Mckendrick S.I.R model; they consider an incidence rate of the form $g(x) y$, where $g$ a function on $[0,\infty)$ that satisfies 

\begin{enumerate}
\item $g(0)=0,$
\item $g(x)\ge 0$ for all $x\in[0,\infty),$
\item There exists $k_g\in\mathbb{R}^{+}$ such that $g(x)\leq k_g,$ for all $x\in[0,\infty).$
\item $g\in C^{1}([0,\infty);\mathbb{R})$, $g^{\prime}(0)>0$, and
\item $g(x)\le g'(0)x$ for all $x\in[0,\infty).$
\end{enumerate}

\noindent The function $g$ introduces potential psychological effects on the population.
This function indicates that in the presence of many infected individuals, the population tends to decrease the number of contacts per unit time. These ideas appeared after  a study on the spread of cholera in 1973. Indeed, in this work, we study in a deterministic and stochastic way the following model

 \begin{equation}\label{sire11}
\noindent 
\begin{aligned}
dy&=( -g(x)y+\gamma x +\mu(N- y))dt,\\
dx&=(g(x)y -(\gamma+\mu)x\Big)dt , t\ge 0. 
\end{aligned}
\end{equation}
  Since $x+y=N$, the equations in (\ref{sire11})  are coupled. Thus, it suffices to analyze one equation. For instance,  

\begin{equation}\label{sis}
    dx=(g(x)(N-x) -(\gamma+\mu)x)dt,\ \ t\ge 0.
\end{equation}

If we assumed that $g(x):=\beta h(x),$ where $h$ satisfies properties $1-5$ and $\beta$ is the coefficient of illness transmission, then the above equation adopts the form

\begin{equation}
\label{sisd}
    dx=(\beta h(x)(N-x) -(\gamma+\mu)x)dt,\ \ t\ge 0.
\end{equation}

In this work, we consider the following random perturbation $\widehat\beta dt=\beta dt+\sigma  dB(t) $, where $\sigma$ is the perturbation rate and $B(t)$ is a one-dimensional Brownian motion. The stochastic differential equation then takes the form 

\begin{equation}\label{siss}
    dx=(\beta h(x)(N-x) -(\gamma+\mu)x)dt+\sigma h(x)(N-x)dB(t),\ \ t\ge 0.
\end{equation}

These kinds of random perturbations are described in \cite{Gray},\cite{TORNATORE2005111}. Particularly in \cite{Xiao}, two types of monotone and non-monotone functions are described. In \cite{Lah}, they deal with the monotone case. 
The function $h(x)$ measures the inhibition effect on susceptible individuals' behavioral changes.

In this paper, we provide sufficient conditions for positivity, stability, extinction and persistence for the solutions of (\ref{siss}). Finally, we show simulations for two specific cases of functions with different behaviors, one of which is monotonic and the other vanishes at infinity.
\section{Main Results}
\label{sec:main}

One of the main properties in stochastic epidemiological models is that the number of individuals of each species remains positive over time. Due to the nature of our model, we also need to ensure that no species exceeds the total population size at any point. 

For the sake of completeness we give the definition of invariant set for stochastic differential equations. Indeed, we say that a set $\mathcal{O}$ is an invariant set for a stochastic differential equation if it has a solution with initial condition in $\mathcal{O}$ which a.s. remains in $\mathcal{O}$ for any $t \geq 0$.

\begin{theorem}\label{Invariance1} The set $\mathcal{O} := (0,N)$ is an invariant set for the stochastic differential equation (\ref{siss}).
\end{theorem}

In the result below we give sufficient conditions such that  if the total population size is not sufficiently large then with high probability the infected population become extincts. 

\begin{theorem}\label{Stability}
Assume that
\begin{equation}\label{condition}
h'(0)<\dfrac{(\gamma+\mu)-\frac{1}{2}\sigma^{2} (h'(0))^{2}N^{2}}{\beta N}.
\end{equation}
Then $x=0$ is an equilibrium point which is asymptotically stable in probability.
\end{theorem}

Nest theorem provides sufficient conditions for global  extinction (a.s.); i.e. there is extinction for any initial condition inside $(0,N)$.

\begin{theorem}\label{Extinction}
If $ h'(0)<\dfrac{\gamma+\mu}{(\beta+\sigma^{2}k_{h})N}$, then for any given initial condition $x(0) \in(0,N)$ the solution of the SDE (\ref{siss}) satisfies
 \begin{equation} \nonumber \displaystyle\mathbb{P}\Big(\displaystyle\lim_{t\to\infty}x(t)=0\Big)=1.
 \end{equation}
\end{theorem}

In the next theorems sufficient conditions under which there is persistence of the disease (a.s.)In the theorem below, we .

\begin{theorem}\label{Spersistence}
Let $h$ be as in (\ref{siss}) and set $\varphi(x)=\dfrac{h(x)}{x}(N-x)$. Assume that
\begin{enumerate}
    \item $\dfrac{\beta }{\sigma^{2}N}<h'(0)<\dfrac{\beta + \sqrt{\beta^2-2\sigma^2(\gamma+\mu)}}{\sigma^{2}N}$,
    \item $\beta^2> 2\sigma^2(\gamma+\mu)$ and
    \item $\varphi'(x)<0,$ for any $x\in (0,N).$ 
\end{enumerate}
Then, for any given initial condition $x(0) \in (0,N)$ of (\ref{siss}),
\begin{enumerate} [(i)] 
    \item \label{lims}
   $\displaystyle\limsup_{t\to\infty}x(t)\geq \xi,  \ \ \mathbb{P}$-a.s.
     \item \label{limi}
   $\displaystyle\liminf_{t\to\infty}x(t)\leq \xi,   \ \ \mathbb{P}$-a.s.
\end{enumerate}
where $\xi\in (0,N)$ and $\varphi(\xi)=\dfrac{\beta - \sqrt{\beta^2-2\sigma^2(\gamma+\mu)}}{\sigma^2}$.
\end{theorem}

Finally, in the following result and under the same conditions as in the previous theorem, we guarantee the existence of only one stationary measure.

\begin{theorem}\label{stationarymeasure}
Let $h$ be as in (\ref{siss}) set $\varphi(x)=\dfrac{h(x)}{x}(N-x)$. Assume that
\begin{enumerate}
\item $\dfrac{\beta }{\sigma^{2}N}<h'(0)<\dfrac{\beta + \sqrt{\beta^2-2\sigma^2(\gamma+\mu)}}{\sigma^{2}N}$, and
\item $\beta^2> 2\sigma^2(\gamma+\mu)$. 
\item  $\varphi'(x)<0,$ for any $x\in (0,N).$ 
\end{enumerate}
Then (\ref{siss}) admits only one stationary measure.
\end{theorem}

\section{Proofs}

Before embarking into the proofs we observe, for the sake of completeness, that the stochastic differential equation (\ref{siss}) has an associated operator called the Lyapunov operator, and it is given by

\begin{equation}
  L_s=\frac{\partial}{\partial t}+\left[\beta h(x)(N-x)-(\gamma+\mu)x\right]\frac{\partial}{\partial x}+0.5 \sigma^2h^2(x)(N-x)^2\frac{\partial}{\partial x^2}.
\end{equation}

\begin{proof} of Theorem \ref{Invariance1}\\

We begin by observing that the functions $b, \sigma : \mathcal{O} \to \mathbb{R}$ given by $b(x)=\beta h(x)(N-x)-(\gamma+\mu)x$ and $\sigma(x)=\sigma h(x)(N-x)$ and the function $h$ (see equation (\ref{siss})) are  locally Lipschitz in $\mathcal{O}$. 

Let $V:\mathcal{O}\to \R$ be given by
\begin{equation}\nonumber
V(x)=\dfrac{1}{x}+\dfrac{1}{N-x}
\end{equation}
and note that $V \in C^{2}( \mathcal{O};[0,+\infty))$. For each $n\in \N$ set $\mathcal{O}_n:=(\frac{1}{n}, N-\frac{1}{n})$. Then $\mathcal{O}_{n} \uparrow \mathcal{O}, \overline{\mathcal{O}_{n}} \subset \mathcal{O}$ and $V_{n}:= \displaystyle\inf_{x\notin \mathcal{O}_{n}}{V(x)} \to +\infty$ as $n\to +\infty$.

Now, for $x\in \mathcal{O},$ we get
 \begin{equation*}
  \begin{aligned}
    \ \ L_{s}V(x)&=b(x)V'(x)+\frac{1}{2}\sigma^{2}(I)V''(x)\\
     &=\dfrac{-\beta h(x)(N-x)}{x^2}+\dfrac{(\gamma+\mu)}{x}+\dfrac{\beta h(x)}{(N-x)}-\dfrac{(\gamma+\mu)I}{(N-x)}\\
     &+\sigma^2 \dfrac{h^2(x)}{x^2}(N-x)^2\dfrac{1}{x}+
     \sigma^2 h^2(x)\dfrac{1}{(N-x)}\\
     & \le \dfrac{(\gamma+\mu)}{x}+\dfrac{\beta k_{h}}{(N-x)}+\sigma^2 (h'(0))^{2}N^{2}\dfrac{1}{x}+\sigma^2 k_{h}^{2}\dfrac{1}{(N-x)}\\
 \end{aligned}   
 \end{equation*}
 where $L_s$ is the Lyapunov operator associated to the SDE (\ref{siss})

Set $c:=\max\{\gamma+\mu+\sigma^2 (h'(0))^{2}N^{2},\beta k_{h}+\sigma^{2}k_{h}^2\}$. Then
$$L_{s}V(x)\leq c \Big(\dfrac{1}{x}+\dfrac{1}{N-x}\Big)=c V(x).$$
It follows from Theorem $3.5$, page $75$ in \cite{khasminskii} 
that for any initial condition the SDE (\ref{siss}) admits a unique solution  $I$  such that 
\begin{equation} \nonumber
\mathbb{P}\Big( x(t)\in (0,N),\forall t\geq 0\Big)=1.
\end{equation} 
\end{proof}

\begin{proof} of Theorem \ref{Stability}.

Since $h(0)=0$ we have that $x=0$ is an equilibrium point for (\ref{siss}). We call this point the disease-free equilibrium point. Now we proceed to prove that $x=0$ is asymptotically stable in probability; i.e. we will show that with high probability the solution to equation (\ref{siss}) and initial condition $x_0=0$ is arbitrary close to $0$.

For $x \in [0,N)$ let $b(x):=\beta h(x)(N-x)-(\gamma+\mu)x$ and $\sigma(x):=\sigma h(x)(N-x)$.

Then
\begin{equation}\nonumber
b'(0)=\beta h'(0)N-\beta h(0)-(\gamma+\mu)=\beta h'(0)N-(\gamma+\mu)
\end{equation}
and 
\begin{equation}\nonumber 
\sigma'(0)=\sigma h'(0)N-\sigma h(0)=\sigma h'(0)N.
\end{equation}
Now we can prove that  $x=0$ is an asymptotically stable equilibrium point in probability for
\begin{equation}\label{sissl}
    dx=b'(0)xdt+\sigma'(0)xdB(t),\ \ t\ge 0
\end{equation}

\noindent Since the function $V:[0,N)\to \R$ defined by $V(x):=x^{2}$ is positive defined at $x=0$ and $V \in C^{2}([0,N);\mathbb{R})$ then, for any $x\in [0,N)$, we have that
\begin{equation*}
  \begin{aligned}
     L_{s}V(x)=b'(0)xV'(x)+\frac{1}{2}(\sigma'(0)x)^{2}V''(x)
     &= 2b'(0)x^{2}+(\sigma'(0))^{2}x^2 \\
     &=(2b'(0)+(\sigma'(0))^{2})x^2\\
     &= C x^2
 \end{aligned}
\end{equation*}
where $C:=2b'(0)+(\sigma'(0))^{2}=2\beta h'(0)N-2(\gamma+\mu)+\sigma^{2} (h'(0))^{2}N^{2}<0.$ Therefore, $L_{s}V$ is negative defined in $[0,N)$. Thus, $x=0$ is a stable asymptotically equilibrium point in probability for (\ref{sissl}).\\

\noindent Now, let $\epsilon>0$ be fixed. Take $\delta:=\min\{\dfrac{\epsilon}{(NL(\beta+\sigma)+2h'(0) )},h\}$ and $x \in (0,\delta)$. A straightforward computation gives
\begin{equation*}
    \begin{aligned}
 \end{aligned}
     |b(x)-b'(0)x|+|\sigma(x)-\sigma'(0)x|\leq \epsilon x.
\end{equation*}

 \noindent It follows from the Lyapunov linearization theorem (See Theorem $7.1$, page $114$ in \cite{khasminskii}) that $x=0$ is an asymptotically stable  equilibrium point in probability  for (\ref{siss}).
\end{proof}

\begin{proof} of Theorem \ref{Extinction}\\
We begin by pointing out that in order to prove Theorem (\ref{Extinction}) it suffices to show that
 \begin{equation} \nonumber 
 \displaystyle\limsup_{t\to\infty}\frac{\ln(x(t))}{t}<0,\ \ \mathbb{P}-a.s.
 \end{equation}
 We know from the proof of the Theorem that there exists a (unique) solution $I$ to (\ref{siss}) defined in $[0,\infty)$ such that $x(t)\in (0,N)$ a.s.for any $t\in[0,\infty)$. 
 
 Now we introduce the function $V:(0,N)\to\mathbb{R} \in C^{2}((0,N) \mathbb{R})$ defined by $V(x)=\ln(x)$. It follows from It\^o formula that
\begin{equation} \nonumber V(x(t))-V(x(0))=\displaystyle\int_{0}^t L_{s}V(x(s) )ds+\displaystyle\int_{0}^t V'(x(s))\sigma(x(s))dB(s) \ \ {\text{a.s}}
 \end{equation}
 Therefore,
 \begin{equation*}
    \begin{aligned}
   \ln(x(t))&=\ln(x_0)+\displaystyle\int_{0}^t L_{s}V(x(s))ds+\displaystyle\int_{0}^t V'(x(s))\sigma(x(s))dB(s)\\
    &=\ln(x_0)+\int_{0}^t\Big(\frac{\beta h(x(s))}{x(s)}(N-x(s)) -(\gamma+\mu)\Big) \\
    & -\frac{1}{2}\frac{\sigma^{2}h^{2}(x(s))}{x^{2}(s)}(N-x(s))^2 \ ds  +\int_{0}^t\frac{\sigma h(x(s))}{x(s)}(N-x(s))dB(s)\\
   &\leq\ln(x_0)+(\beta h'(0)N-(\gamma+\mu)+\sigma^{2} h'(0)k_{h}N)t \\
   & +\int_{0}^t\frac{\sigma h(x(s))}{x(s)}(N-x(s))dB(s),\ \ {\text{a.s}.} 
 \end{aligned} 
 \end{equation*}

By the Law of large numbers for martingales we have
\begin{equation} \nonumber \displaystyle\lim_{t\to\infty}\frac{1}{t}\int_{0}^t\frac{\sigma h(x(s))}{x(s)}(N-x(s))dB(s)=0\ \ {\text{a.s}.}
\end{equation}
Thus, 
\begin{equation} \nonumber 
\limsup_{t\to\infty}\frac{\ln(x(t))}{t}\leq\beta h'(0)N-(\gamma+\mu)+\sigma^{2} h'(0)k_{h}N
<0,\ \ {\text{a.s}.}
\end{equation}
Thence, $\displaystyle\lim_{t\to\infty}x(t)=0,$ a.s. This finishes the proof.
\end{proof}


\noindent {\bf{Proof of Theorem \ref{Spersistence}}} Olhar o tex nessa linha.

\begin{proof} 
    Now we will study the persistence regime for the SDE (\ref{siss}). Just consider $\mathcal{O}:=(0,N)$ and let $V: \mathcal{O}\to \mathbb{R}$ be defined by



\[
V(x) := \dfrac{\xi}{N-\xi}\left(\dfrac{N}{x} - 1\right) - 1 - \ln \left[\dfrac{\xi}{N-\xi} \left(\dfrac{N}{x} - 1\right)\right]
\]

Namely, we will prove that for any $t_0 > 0$, $I(t)$ visits infinitely many times the level $\xi$ for $t\geq t_0$ for some $\xi \in (0,N)$ with probability one. 

Let $V(x)=\ln(x)$. Then the associated Lyapunov operator to the SDE (\ref{siss}) is given by

\begin{eqnarray}\nonumber
L_{s}V(x)=-\dfrac{1}{2}\sigma^{2}\left[\varphi(x)-\frac{\beta + \sqrt{\beta^2-2\sigma^2(\gamma+\mu)}}{\sigma^{2}}\right]\left[\varphi(x)-\frac{\beta - \sqrt{\beta^2-2\sigma^2(\gamma+\mu)}}{\sigma^{2}}\right] \\
\end{eqnarray}
where $\varphi(x)=\dfrac{h(x)}{x}(N-x)$. A straightforward but tedious computation yields the existence of $\xi\in(0,N)$ such that $L_sV(\xi)=0$ whenever $\beta^2\ge 2\sigma^2(\gamma+\mu)$ and
$\dfrac{\beta }{\sigma^{2}N}<h'(0)<\dfrac{\beta + \sqrt{\beta^2-2\sigma^2(\gamma+\mu)}}{\sigma^{2}N}$. 

Now we list some properties of the Lyapunov operator $L_sV$ which will be useful in study of persistence. Indeed, if $\varphi'(x)<0$ for any $x\in (0,N)$ then there exist $m\in(0,\xi)$ such that

\begin{enumerate}
    \item  $L_{s}V$ is increasing on $(0,m)$ and $L_{s}V(x)>0$ for any $x\in(0,m).$
     \item $L_{s}V$  is decreasing on $(m,\xi)$ and $L_{s}V(x)>0$ for any $x\in(m,\xi).$
     \item  $L_{s}V$ is decreasing on $(\xi,N)$ and $L_{s}V(x)<0$ for any $x\in(\xi,N).$
\end{enumerate}
where $L_sV(m)$ is the maximum value of $L_sV$ on $(0,N).$

\begin{figure}[htb]
\centering
\includegraphics[width=0.5\textwidth]{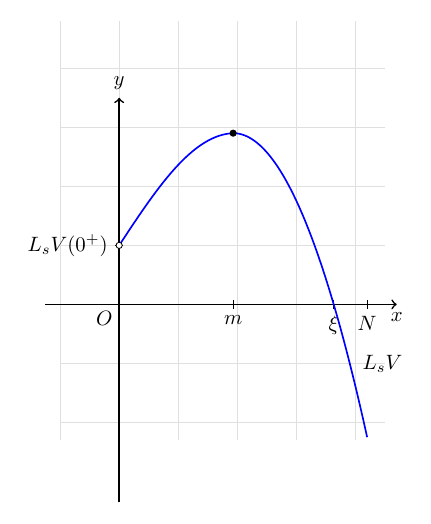}
\caption{Lyapunov operator $L_sV$ for $V(x)=\ln(x)$}
\label{fig:G2}
\end{figure}

In the proof of this theorem we follow the approach developed in \cite{Gray}. Suppose that (\ref{lims}) does not hold. Then there exists $\epsilon\in(0,1)$ such that

\begin{equation} \nonumber
\mathbb{P}\left\{\omega\in\Omega:\displaystyle\limsup_{t\to\infty}x(t,\omega)< \xi\right\}>\epsilon .
\end{equation}

\noindent
Consider a sequence of sets $\left(A_n\right)_{n \geq 1}$ defined by
 
\begin{equation} \nonumber
A_n:=\{\omega\in\Omega:\displaystyle\limsup_{t\to\infty}x(t,\omega)\le \xi-2/n\}.
\end{equation}
 
It is easy to see that $A_n\uparrow\left\{\omega\in\Omega:\displaystyle\limsup_{t\to\infty}x(t,\omega)< \xi\right\}$. Then, there exists $p\in\mathbb{N}$ such that $\mathbb{P}(A_{p})>\epsilon,$ $m<\xi-1/p$ and $L_{s}V(0^{+})>L_{s}V(\xi-1/p)$.
 
Set
\begin{equation} \nonumber
\tilde{\Omega}:=\{\omega\in\Omega:\displaystyle\lim_{t\to\infty}\frac{1}{t}\int_{0}^t\frac{\sigma h(x(s,\omega))}{x(s,\omega)}(N-I(s,\omega))dB(s,\omega)=0\}. 
\end{equation}
It follows from the law of large numbers for continuous time Martingales that $\mathbb{P}(\tilde{\Omega})=1$. Take $B=\tilde{\Omega}\cap A_{p}$ and $\omega\in B$. Then, there exists $T(\omega)>0$ such that

\begin{equation}\label{l1}
 x(t,\omega)<\xi-1/p \hspace{.2cm}\text{for all}\hspace{.2cm}t\geq T(\omega).  
\end{equation}

Since $L_{s}V(0^{+})>L_{s}V(\xi-1/p)$, it follows from the properties of the Lyapunov operator associated to $\ln(x)$ discussed immediately after the statement of Theorem \ref{Spersistence} and from (\ref{l1}) that
\begin{equation} \nonumber
L_{s}V(x(t,\omega))\geq L_{s}V(\xi-1/p).
\end{equation}
By the It\^o's formula we obtain
\begin{equation*}
    \begin{aligned}
   \ln(x(t,\omega))&=\ln(x_0)+\displaystyle\int_{0}^t L_{s}V(x(s,\omega))ds+\int_{0}^t\frac{\sigma h(x(s,\omega))}{x(s,\omega)}(N-x(s,\omega))dB(s,\omega)\\
    &=\ln(x_0)+\displaystyle\int_{0}^{T(\omega)} L_{s}V(x(s,\omega))ds+\int_{{T(\omega)}}^{t} L_{s}V(x(s,\omega))ds\\
    &+\int_{0}^t\frac{\sigma h(x(s,\omega))}{x(s,\omega)}(N-x(s,\omega))dB(s,\omega)\\
     &\geq\ln(x_0)+\displaystyle\int_{0}^{T(\omega)} L_{s}V(x(s,\omega))ds+L_{s}V(\xi-1/p)(t-T(\omega))\\
     &+\int_{0}^t\frac{\sigma h(x(s,\omega))}{x(s,\omega)}(N-x(s,\omega))dB(s,\omega)
 \end{aligned} 
 \end{equation*}

Dividing by $t$ and taking inferior limit we get
\begin{equation} \nonumber
  \liminf_{t\to\infty}\frac{\ln(x(t,\omega))}{t}\geq L_{s}V(\xi-1/p)>0.
\end{equation}
Therefore, $\lim_{t\to\infty}x(t,\omega)=\infty$ which contradicts (\ref{l1}). 

Now, suppose that (\ref{limi}) does not hold. Then, there exists $\epsilon\in(0,1)$ such that 
\begin{equation} \nonumber 
\mathbb{P}\{\omega\in\Omega:\displaystyle\liminf_{t\to\infty}x(t,\omega)> \xi\}>\epsilon.
\end{equation}
Consider the sequence of sets $\left(C_n\right)_{n \geq 1}$ defined by  
\begin{equation*}
C_n:=\{\omega\in\Omega:\displaystyle\liminf_{t\to\infty}x(t,\omega)\geq \xi+2/n\}.
\end{equation*}
It is easy to see that $C_n\downarrow\{\omega\in\Omega:\displaystyle\liminf_{t\to\infty}x(t,\omega)> \xi\}$. Therefore,
there exists $p\in\mathbb{N}$ such that $\mathbb{P}(C_p)>\epsilon$. Let and $\omega\in C=\tilde{\Omega}\cap C_p$. Then, there exists $T(\omega)>0$ such that for any $t\geq T(\omega)$ the following inequality holds
\begin{equation}\label{l2}
x(t,\omega)>\xi+1/p  .
\end{equation}
It follows, again, from the properties of the Lya-
punov operator associated to ln(x) discussed immediately after the statement of Theorem \ref{Spersistence} and from (\ref{l1}) that
\begin{equation} \nonumber
L_{s}V(x(t,\omega))\leq L_{s}V(\xi+1/p).
\end{equation}

Now, by the It\^o formula, we have
\begin{equation*}
    \begin{aligned}
   \ln(x(t,\omega))&=\ln(x_0)+\displaystyle\int_{0}^t   L_{s}V(x(s,\omega))ds  \\
   &+\int_{0}^t\frac{\sigma h(x(s,\omega))}{x(s,\omega)}(N-x(s,\omega))dB(s,\omega)\\
    &=\ln(x_0)+\displaystyle\int_{0}^{T(\omega)} L_{s}V(x(s,\omega))ds+\int_{{T(\omega)}}^{t} L_{s}V(x(s,\omega))ds\\
    &+\int_{0}^t\frac{\sigma h(x(s,\omega))}{x(s,\omega)}(N-x(s,\omega))dB(s,\omega)\\
    & \leq\ln(x_0)+\displaystyle\int_{0}^{T(\omega)} L_{s}V(x(s,\omega))ds+L_{s}V(\xi+1/p)(t-T(\omega))\\
     & +\int_{0}^t\frac{\sigma h(x(s,\omega))}{x(s,\omega)}(N-x(s,\omega))dB(s,\omega).
 \end{aligned} 
 \end{equation*}
 Then 
 \begin{equation} \nonumber \limsup_{t\to\infty}\frac{\ln(x(t,\omega))}{t}\leq L_{s}V(\xi+1/p)
<0.
\end{equation}
Therefore, $$\lim_{t\to\infty}x(t,\omega)=0$$ which contradicts (\ref{l2}). This finishes the proof.
 
\end{proof}

\begin{proof} of Theorem \ref{stationarymeasure}\\

Before beginning with the proof of the existence and uniqueness of the invariant measure we briefly remind its definition. Let $\left(x(t)\right)_{t \geq 0}$ be a a solution of an SDE with initial condition $x(0)$. Denote by  ${\mathbb{P}}_t^\mu$ the law of $x(t)$ and assume that $x(0)\sim\mu$. As usual we say that the measure $\mu$ is invariant for $x(t)$ if $
{\mathbb{P}}_t^\mu(\cdot)=\mu(\cdot)$ for any $t > 0$.

In \cite{khasminskii} the authors provide sufficient conditions for the existence of a unique stationary distribution given a stochastic process. In particular, we can found the following lemma: 

\begin{lemma}\label{lemmastationarymeasure}
Let $a, b$ with $a<b$. Set $\tau:=\inf\{t\ge 0: x(t)\in(a,b)\}$. If $\mathbb{E}(\tau) <\infty$ for any  $x_0\in(0,a]\cup[b,N)$, then the stochastic process $\{x(t)\}_{t\geq 0}$ admits a unique stationary measure.
\end{lemma}

Thus, to prove this result it will be sufficient to prove that the solution of (\ref{siss}) satisfies the Lemma \ref{lemmastationarymeasure}. Now we prove Theorem \ref{stationarymeasure}.

Let $\xi \in (0,N)$ be as in the conclusion of Theorem \ref{Spersistence}. Take $a, b \in (0,N)$ such that $a< \xi < b$. It follows from the properties of $L_s V$ listed  before Theorem \ref{Spersistence} that

\begin{equation}\label{eccc1}
   L_sV(x)\ge \min\{L_{s}V(0^{+}),L_{s}V(a)\} \text{for any}\ \ x\in(0,a]
\end{equation}\label{ecc1}
and that
\begin{equation}
L_sV(x)\le L_{s}V(b)\text{for any}\ \ x\in[b,N).
\end{equation}
Consider the stopping time

\begin{equation} \nonumber 
\tau:=\inf\{t\ge 0: x(t)\in(a,b)\}.
\end{equation}
One may assume without loss of generality that the initial condition $x_0\in(0,a)\cup (b,N)$. If $x_0\in(0,a)$, then

\begin{equation}
  \ln(a)\ge \ln(x(t\land\tau))) \hspace{.2cm} \text{for any} \ \ t>0,\mathbb{P}-a.s.
\end{equation}
It follows from Dynkin Lemma and inequality  (\ref{eccc1}) that
\begin{equation} \nonumber
\mathbb{E}[\ln(x(t\land\tau))]\ge \ln(x_0)+\min\{L_{s}V(0^{+}),L_{s}V(a)\}\mathbb{E}(t\land\tau).
\end{equation}
Therefore
\begin{equation} \nonumber
\ln(a)\ge \ln(x_0)+\min\{L_{s}V(0^{+}),L_{s}V(a)\}\mathbb{E}(t\land\tau).
\end{equation}
Then, by passing to the limit when $t$ goes to infinity, we have
\begin{equation}\label{imp1}
    \mathbb{E}(\tau)\le \frac{\ln(a/x_0)}{\min\{L_{s}V(0^{+}),L_{s}V(a)\}}
\end{equation}
On the other hand if $x_0\in(b,N)$ then, by a similar argument, we get
\begin{equation} \nonumber 
\ln(b)\le \ln(x_0)-|L_{s}V(b)|\mathbb{E}(t\land\tau).
\end{equation}

Then, by passing to the limit when $t$ goes to infinity, we arrive to
\begin{equation}\label{imp2}
\mathbb{E}(\tau)\le \frac{\ln(N/b)}{|L_{s}V(b)|}.
\end{equation}
Thus we may conclude, by means of (\ref{imp1}) and (\ref{imp2}) that $\mathbb{E}(\tau)<\infty$ for any $x_0\in(0,a]\cup[b,N)$. The desired conclusion follows from Lemma \ref{lemmastationarymeasure}.

\end{proof}

\section{Examples and Simulations}
Mathematically speaking, the psychological effect is responsible for a reduction in the incidence rate $\beta h(x)y$ and in the total number of contacts. For further details on this psychological phenomenon, see \cite{Xiao}.

\subsection{Simulations for $h_1$}
In this subsection, we show through simulations how extinction occurs in the deterministic and stochastic cases for a monotonic function $h_1$. Also, we show through simulations how the phenomenon of stochastic persistence occurs for a monotonic function $h_1$. Finally, since we do not know an explicit formula for the stationary measure ,we exhibit a histogram of the stationary measure for a monotonic function $h_1$.
\subsubsection{Exctintion}

In both deterministic and stochastic cases, we simulate the disease's extinction for a monotone function $h_1$. We observe how the disease disappears after a sufficiently long time.
 
\begin{figure}[H]
 \centering
 \subfloat[Deterministic Extinction. ]{
   \label{f:DetExt}
    \includegraphics[width=0.5\textwidth]{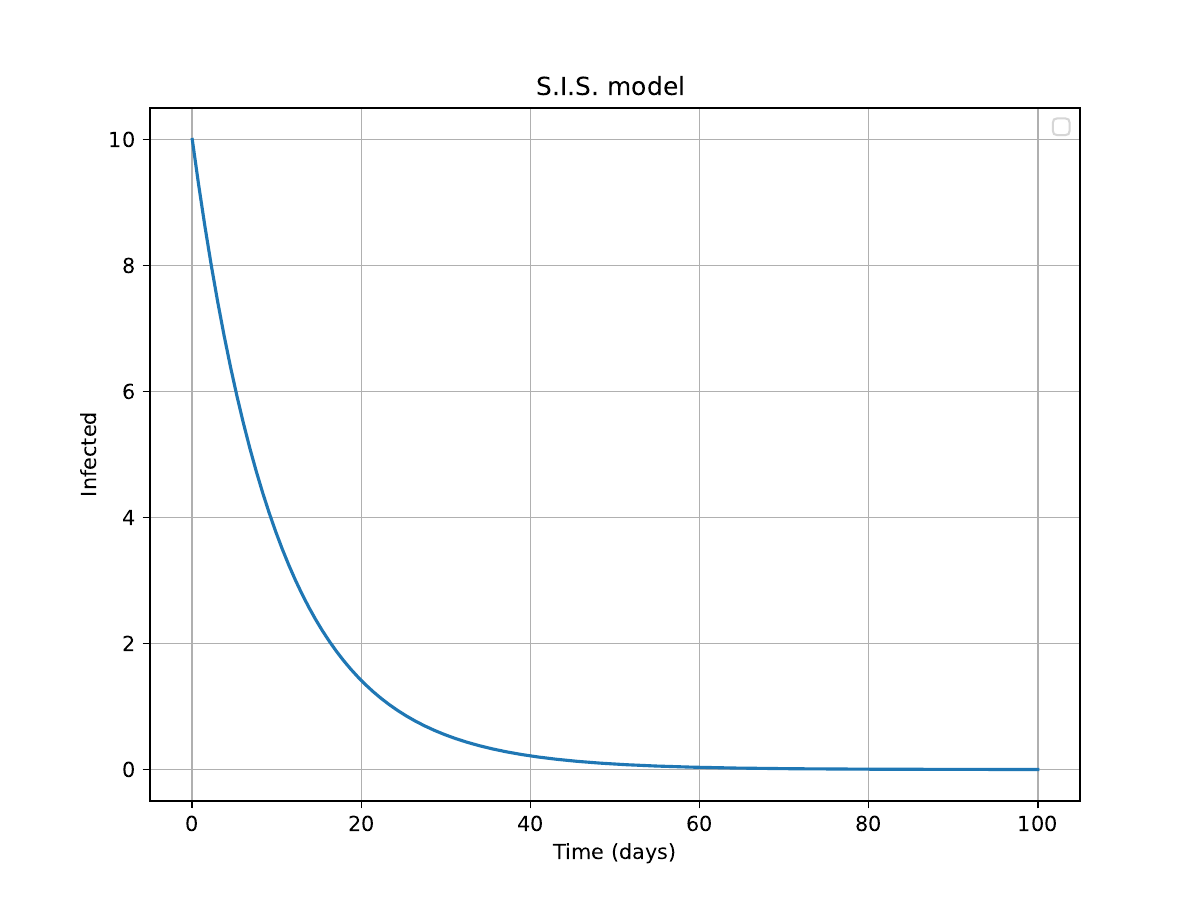}}
  \subfloat[Stochastic Extinction.
]{
   \label{f:stoext}
    \includegraphics[width=0.5\textwidth]{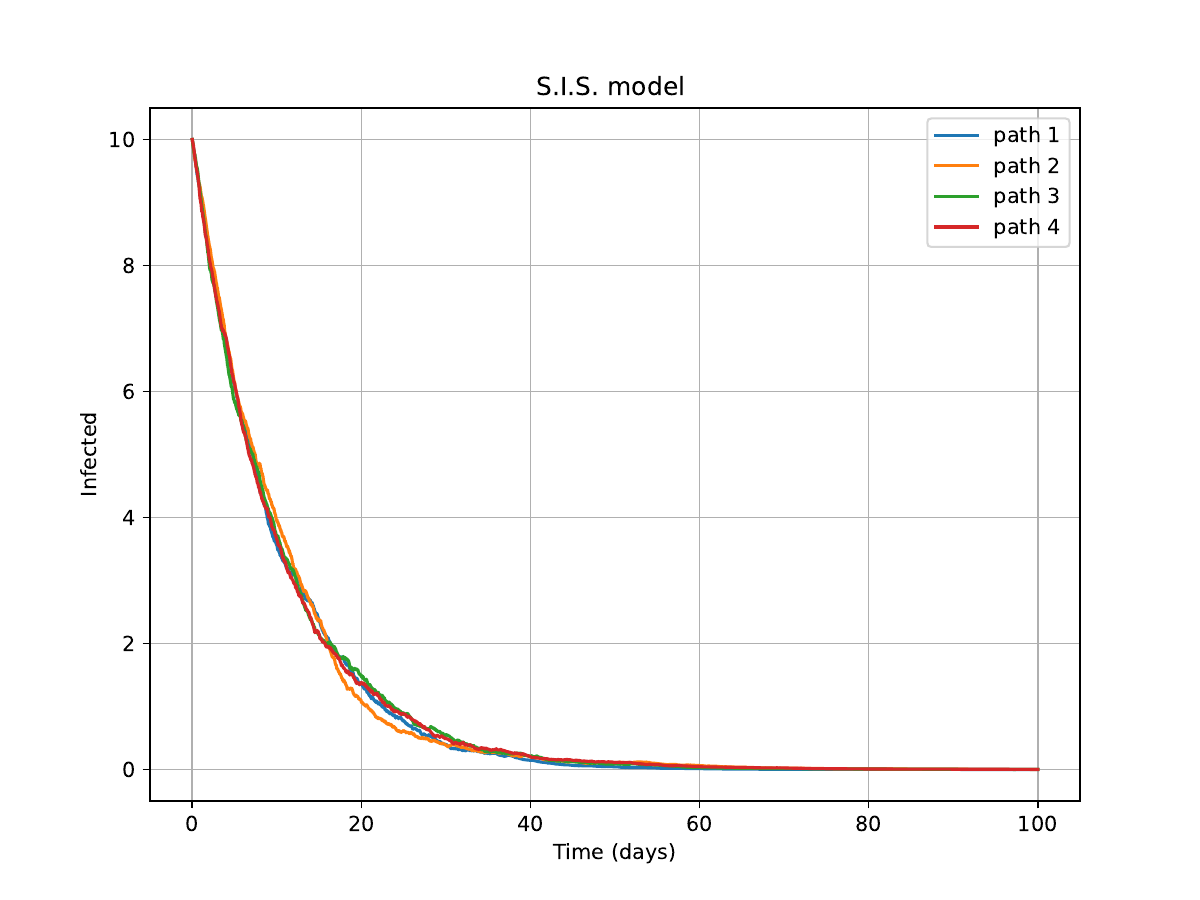}}  
    \caption{Stochastic and deterministic extinction}
 \label{f:extinction}
\end{figure}

 For the deterministic case, in Figure (\ref{f:DetExt}), we take $ \beta= 0.00001, \gamma=0.1, \mu=0.0001, N=1000$, and $x_0=10$. In the stochastic case, in Figure (\ref{f:stoext}), the perturbation was introduced with $\sigma=0.00009$. The rates were $\beta=0.00001, \gamma=0.1, \mu=0.0001$, $N=1000$, and $x_0=10$.

\subsubsection{Persistence}
\begin{figure}[H]
 \centering
  \subfloat[Stochastic persistence $x_0=1$
]{
   \label{f:Assta}
    \includegraphics[width=0.5\textwidth]{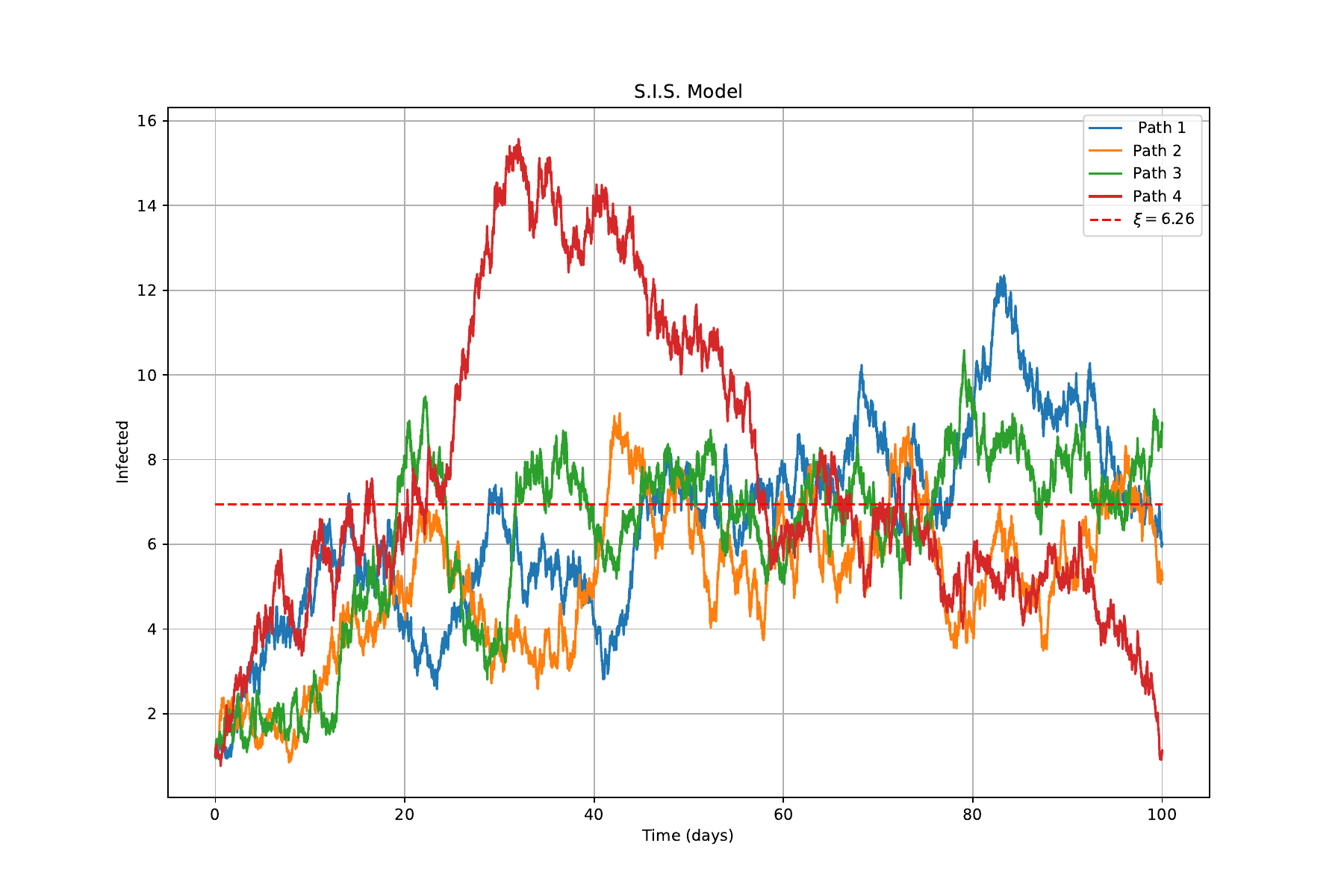}}
  \subfloat[Stochastic persistence $x_0=100$]{
   \label{f:stoper}
    \includegraphics[width=0.5\textwidth]{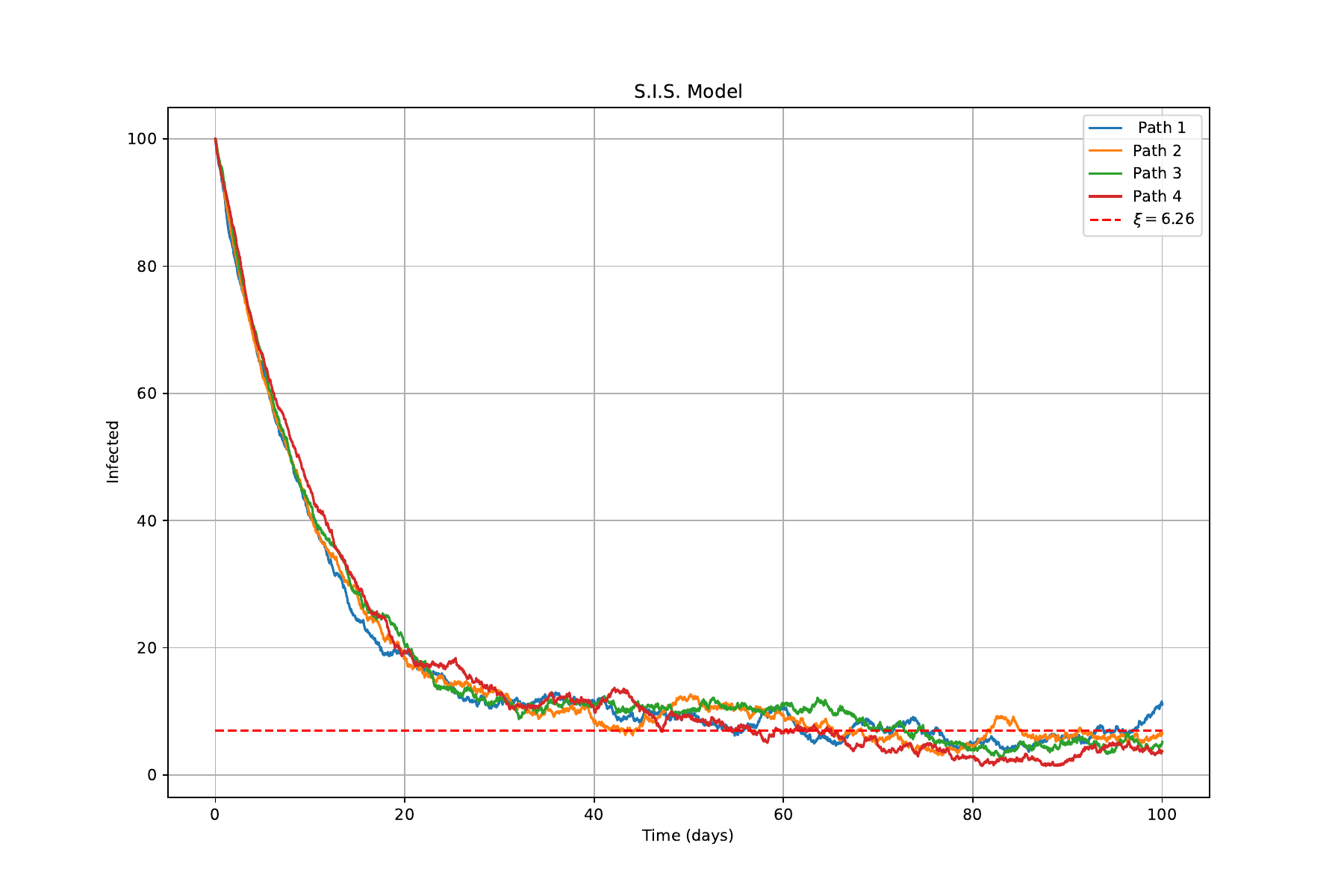}}
    \caption{Stochastic persistence}
 \label{f:extinction}
\end{figure}

Now, we simulate the stochastic persistence of a disease for the function $h_1$. We observe how the disease is recurrent at an endemic level $\xi$. We present two perspectives for different initial conditions.
We show a simulation for $\sigma=0.001$, $\beta= 0.0008, \gamma=0.1, \mu=0.0001, N=1000$.\\
In the figure (\ref{f:Assta})we have $x_0=1$ and in the figure (\ref{f:stoper}) we have $x_0=100.$

\begin{figure}[H]
 \centering
  \subfloat[Stationary distribution $h_1$ function
]{
   \label{f:SD}
    \includegraphics[width=0.5\textwidth]{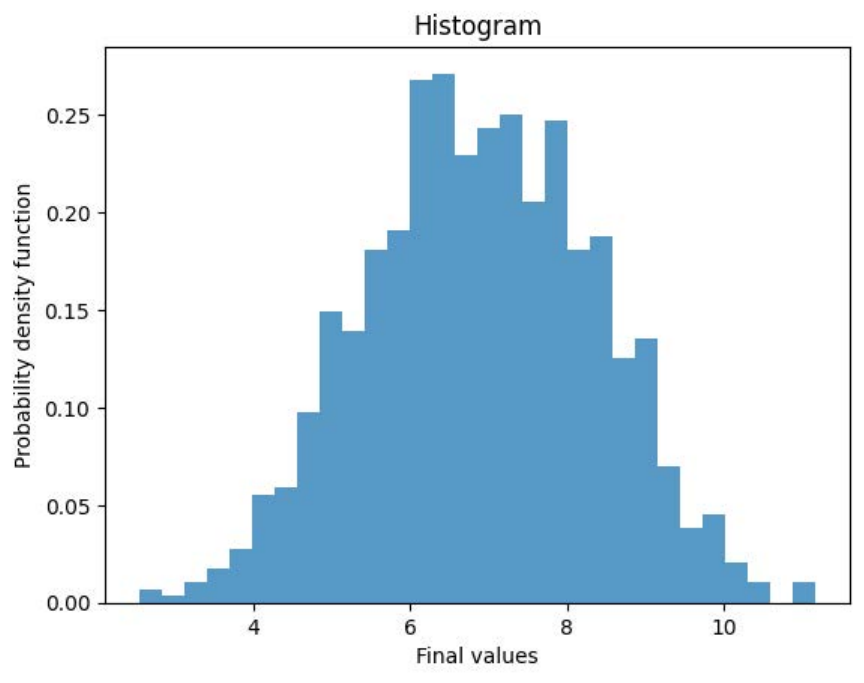}}
    \caption{Stationary distribution}
\end{figure}

\subsection{Simulations for $h_2$}
In this subsection, we show through simulations how extinction occurs in the deterministic and stochastic cases for a non-monotonic function $h_2$. Also, we show through simulations how the phenomenon of stochastic persistence occurs for a non-monotonic function $h_2$. Finally, since we do not know an explicit formula for the stationary measure, we exhibit a histogram of the stationary measure for a non-monotoni function $h_2$.

\subsubsection{Exctintion}
In both deterministic and stochastic cases, we simulate the disease's extinction for a non-monotonic function $h_2$. We observe how the disease disappears after a sufficiently long time.
 
\begin{figure}[H]
 \centering
 \subfloat[Deterministic extinction]{
   \label{f:detext1}
    \includegraphics[width=0.5\textwidth]{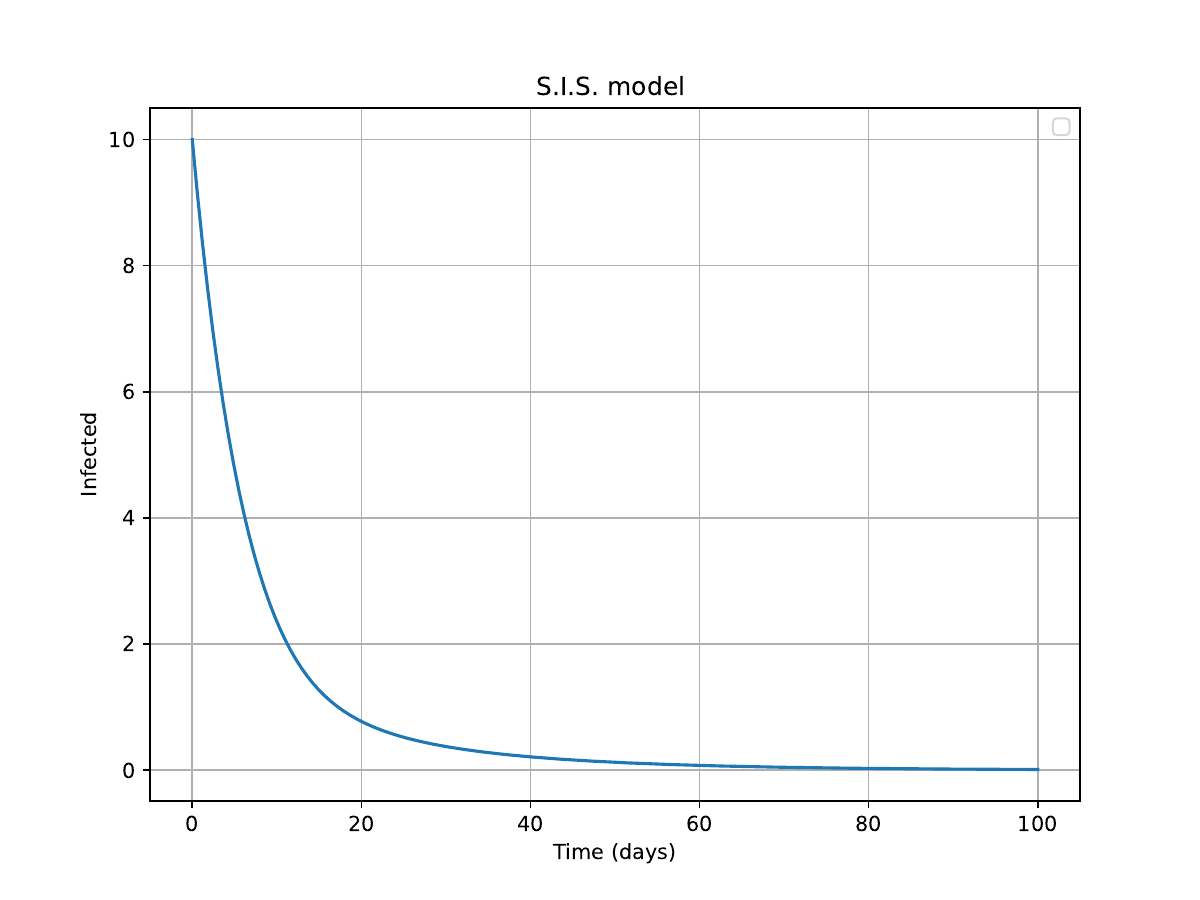}}
  \subfloat[Stochastic extinction
]{
   \label{f:extsto1}
    \includegraphics[width=0.5\textwidth]{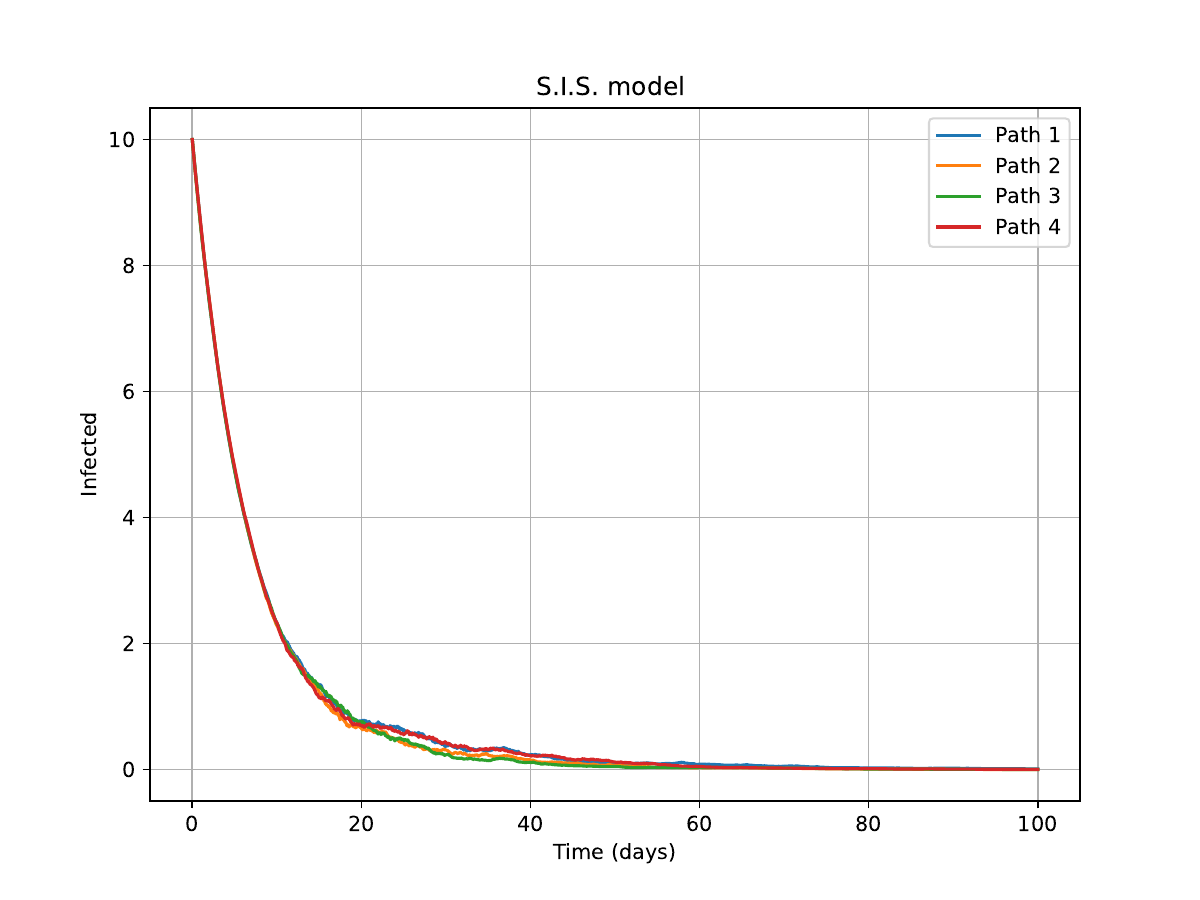}}  
    \caption{Extinction}
 \label{f:extinction}
\end{figure}

In the stochastic case, in Figure (\ref{f:extsto1}), we have $\sigma=0.0009$, $\beta=0.0001, \gamma=0.1, \mu=0.05$, $N=1000$, and $x_0=10$. In the deterministic case, in Figure (\ref{f:detext1}), we have $\beta=0.0001, \gamma=0.1, \mu=0.05$, $N=1000$, and $x_0=10$.

\subsubsection{Persistence}
Now, we simulate the stochastic persistence of a disease for the function $h_2$. We observe how the disease is recurrent at an endemic level $\xi$. We present two perspectives for different initial conditions.
We simulate with $\sigma=0.001$, $\beta= 0.00099, \gamma=0.1, \mu=0.05, N=1000$.

\begin{figure}[H]
 \centering
  \subfloat[stochastic persistence
with $x_0=1$]{
   \label{f:asympsta}
    \includegraphics[width=0.5\textwidth]{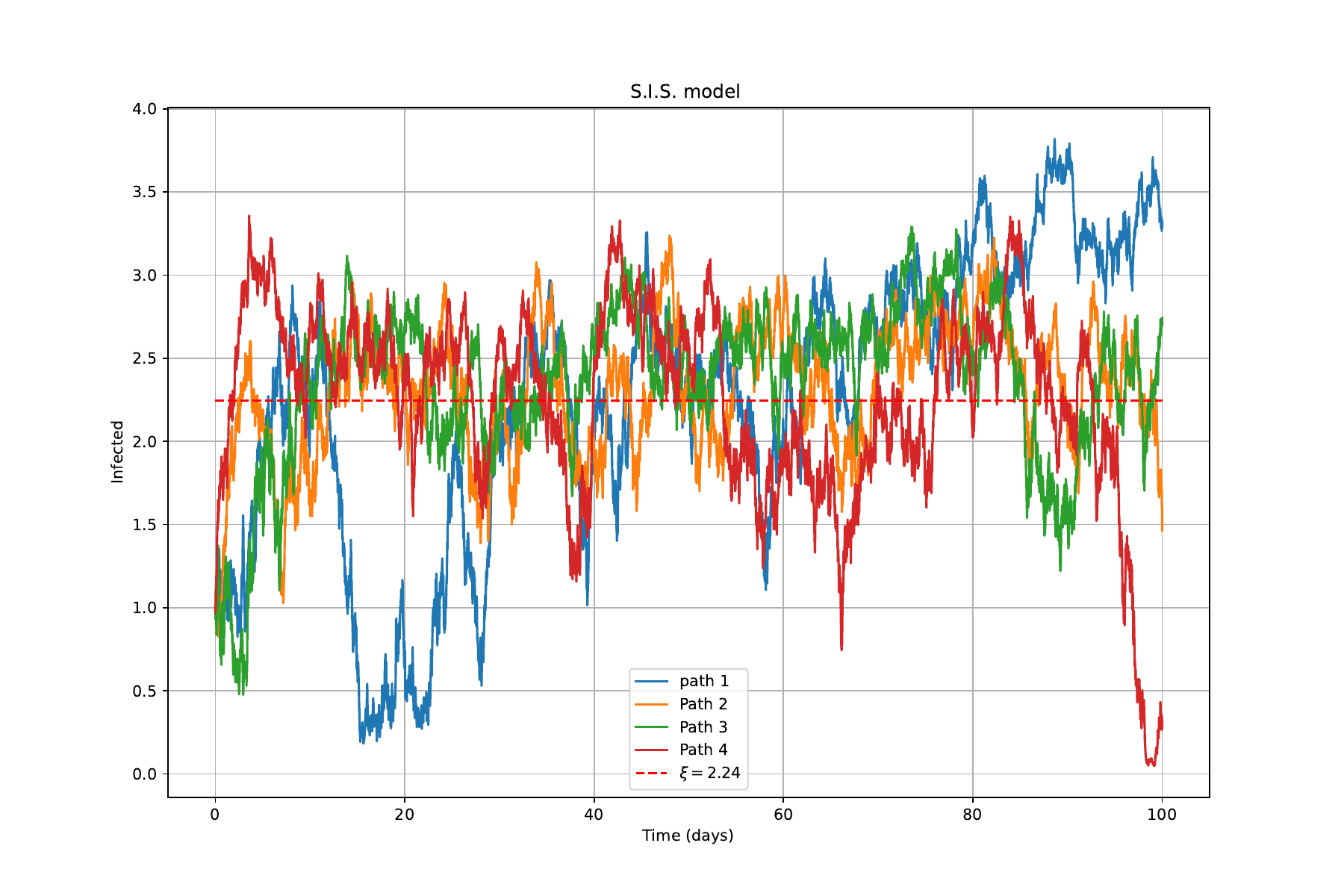}}
  \subfloat[Stochastics persistence with $x_0=100$]{
   \label{f:stoper1}
    \includegraphics[width=0.5\textwidth]{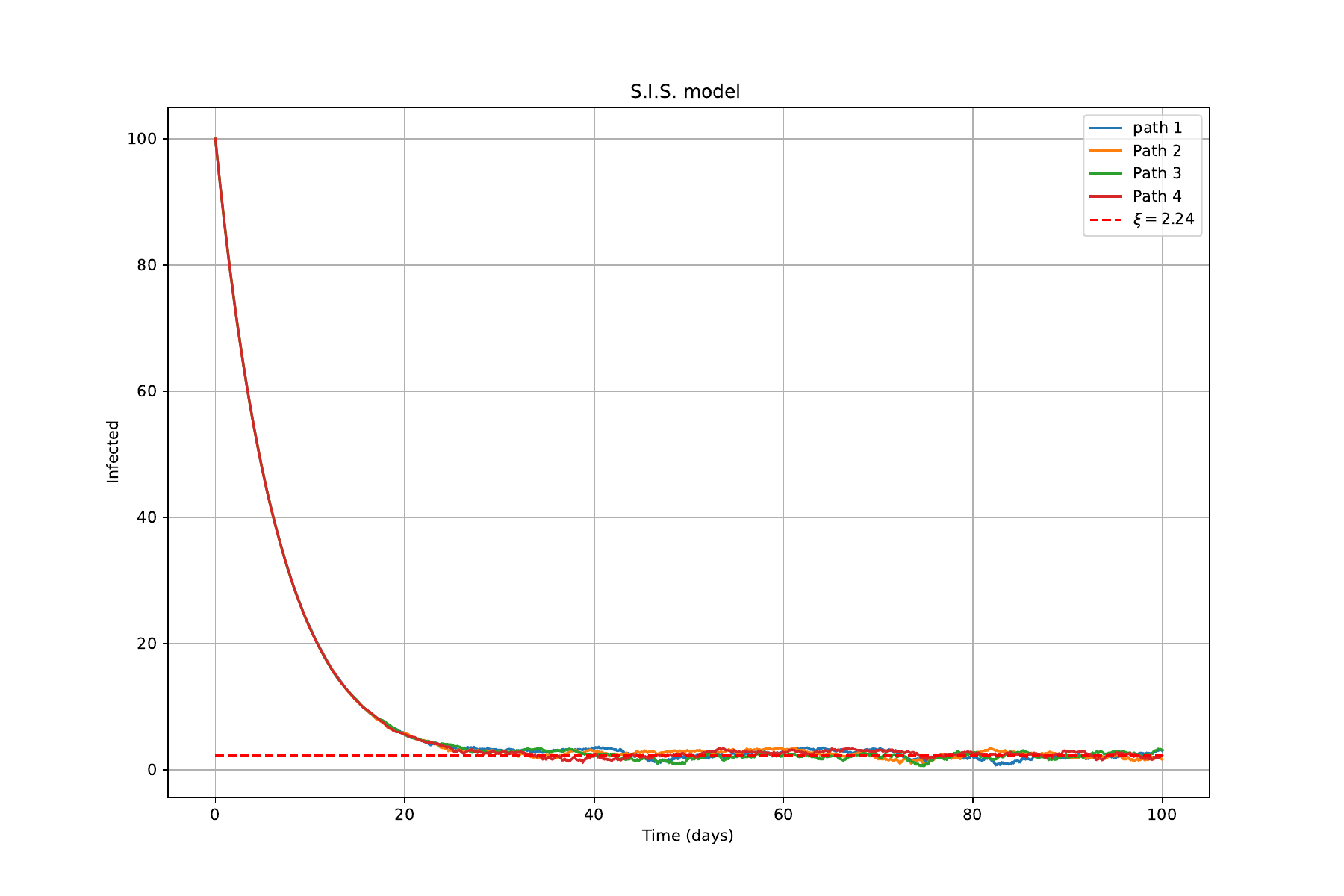}}
    \caption{Stochastic persistence}
 \label{f:extinction}
\end{figure}

In the figure(\ref{f:asympsta}) we have $x_0=1$ and (\ref{f:stoper1}) we have $x_0=100.$

\begin{figure}[H]
 \centering
  \subfloat[Stationary distribution $h_2$ function
]{
   \label{f:SD}
    \includegraphics[width=0.5\textwidth]{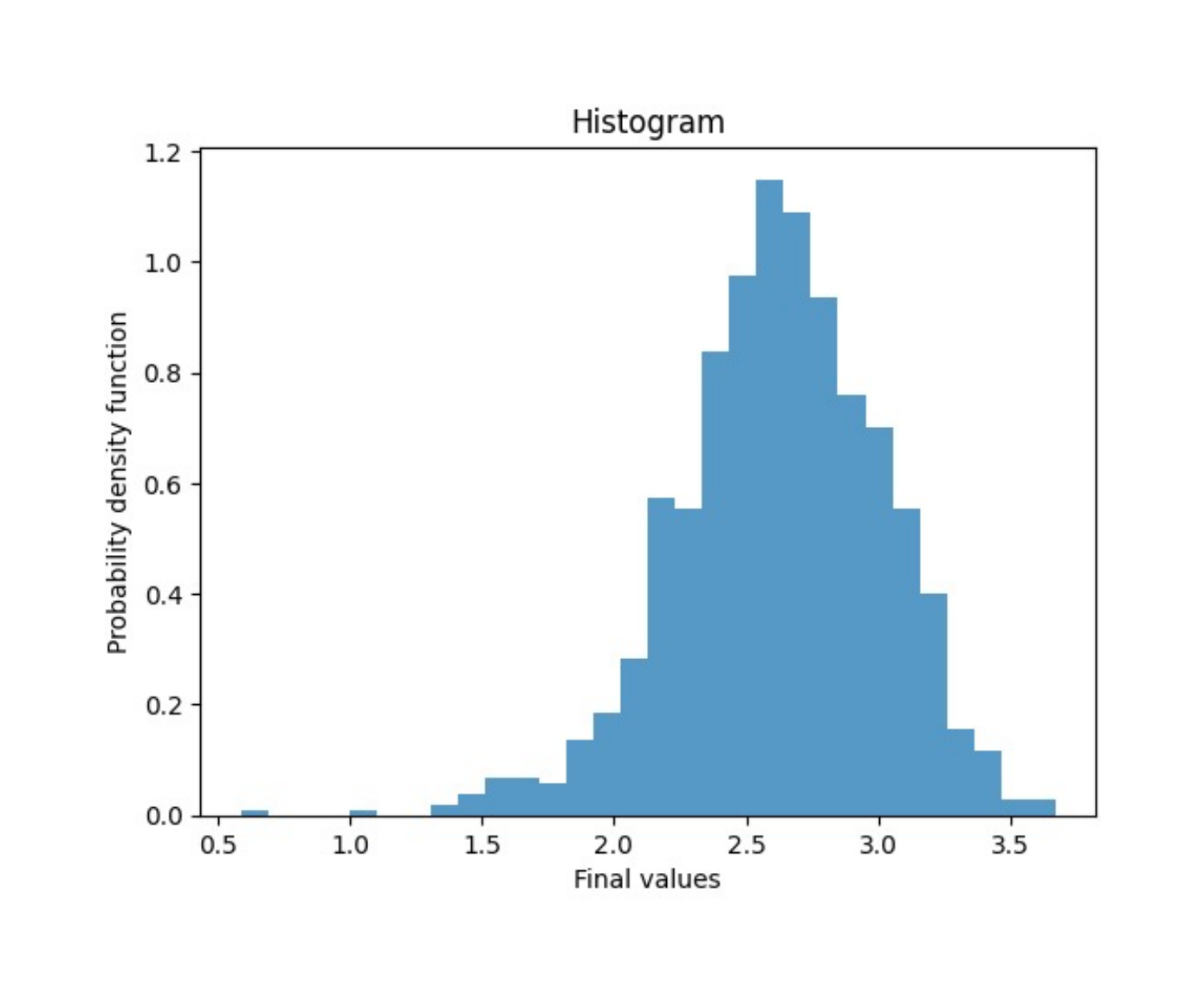}}
    \caption{Stationary distribution}
\end{figure}

\section*{Acknowledgments} 
Research supported by grants 2022/08948-2 and 2017/10555-0 S\~ao Paulo Research Foundation (FAPESP).This research was also supported by grant Universidad de Antioquia - FAPESP 2023-58830. Cristian Coletti thanks UdeA for warm hospitality during his visit in 2023.

\bibliographystyle{plain} 
\bibliography{references} 

 \end{document}